\documentclass[preprint,times,12pt]{article}

\usepackage{amstext}
\usepackage{latexsym}
\usepackage{amssymb}
\usepackage{amsmath}
\usepackage{graphicx}
\usepackage{ifpdf}
\usepackage{color}

\allowdisplaybreaks
\DeclareGraphicsExtensions{.pdf,.png,.jpg,.jpeg,.mps,.eps}

\usepackage[colorlinks=false,linktocpage=true]{hyperref}
\usepackage{hyperref}
\hypersetup{colorlinks=false}

\def\R{\mathbb R}

\def\G{\mathsf{G}}

\def\int{{\rm int\;}}

\def\tr{{\rm tr}}
\def\eps{\varepsilon}
\def\PSL{{\rm PSL}}

\def\X{{\mathfrak X}}
\def\H{\mathbb H}

\def\X{\mathcal X}

\newcommand{\SL}{{\rm SL}}

\newtheorem{theorem}{Theorem}[section]

\newtheorem{definition}[theorem]{Definition}

\newtheorem{lemma}[theorem]{Lemma}

\usepackage{hyperref}

\begin{document}
		\title{A solution to Flinn's conjecture on weakly expansive flows}
	\author{{\sc Huynh Minh Hien} \\[1ex] 
		Faculty of Mathematics and Statistics,\\
		Quy Nhon University,\\
		 170 An Duong Vuong, Quy Nhon, Vietnam;\\[1ex] 
		e-mail: huynhminhhien@qnu.edu.vn  
	}\date{}
	\maketitle 
	\begin{abstract} In L.W. Flinn's PhD thesis published in 1972, the author conjectured that weakly expansive flows are also expansive flows.
		In this paper we use the horocycle flow
		on compact Riemann surfaces of constant negative curvature to show that Flinn's conjecture is not true.
	\end{abstract}

	
	\maketitle
	
	\section{Introduction}
	
	Expansive flows have been being studied since 1972 by several mathematicians with different varieties of expansiveness \cite{artigue3,artigue,bw,flinn,komuro}. Bowen/Walters \cite{bw} and Flinn \cite {flinn} introduced
	the term {\em expansive flow} by different definitions which are equivalent
	to each other in the case of fixed-point-free flows. In his thesis \cite{flinn}, Flinn also provided
	the notion of {\em weak expansiveness} for flows and
	conjectured that weak expansiveness is equivalent to expansiveness. The Flinn conjecture has not been solved until now.
	  In 1984, Gura \cite{gura} discovered a new kind of expansive flows called {\em separating flows}. 
A flow in a compact metric space is said to be {\em separating} if  two orbits of the flow
	must be identical, provided that they are sufficiently close for the whole time. The author  showed that
	every time change of the horocycle flow on a negatively curved compact surface is positive separating. Such a flow is called {\em positive strong separating}. In 2014, Artigue \cite{artigue} provided the concept of 
	{\em kinematic expansiveness} which implies separation. 
	The author also gave the notions of {\em strong kinematic expansiveness}, {\em geometric separation},...
	as well as examples  to analyze relationships among the concepts.  Roughly speaking,
	a flow is strong kinematic expansive if all its  time changes are kinematic expansive. Note that the terms {\em strong kinematic expansiveness} in the sense of Artigue \cite{artigue} and {\em weak expansiveness} in the sense of Flinn \cite{flinn} are equivalent.  In recent paper, Huynh \cite{huynh}
	shows that the horocyle flow on compact Riemann surfaces of constant negative curvature
	is kinematic expansive but not geometric separating (and hence it is not expansive in the sense of Bowen-Walters).  In this paper we prove that the horocycle flow on  compact Riemann surfaces with 
	constant negative curvature is positive strong kinematic expansive. Consequently, a negative answer to Flinn's conjecture is given.
	
	The paper is organized as follows. In the next section, we construct an equivalent flow of the horocycle flow and background material which is well-known in principle. Then in  Section 3 we recall the definitions of some expansive flows and present the strong kinematic expansiveness for the horocycle flow (Theorem \ref{pne}), which followed by a negative answer to Flinn's conjecture (Theorem \ref{fca}).
	\section{Compact Riemann surfaces of constant negative curvature}
	
	We consider the horocycle flow on compact Riemann surfaces of constant negative curvature. 
	It is well-known that any compact orientable surface 
	with constant negative curvature is isometric to a factor $\Gamma\backslash \H^2=\{\Gamma z, z\in\H^2 \}$, where
	$\Gamma $ is a discrete subgroup of the group $\PSL(2,\R)=\SL(2,\R)/\{\pm E_2 \}$;
	here $\SL(2,\R)$ is the group of all real $2\times 2$ matrices with unity determinant, $E_2$ is the unit matrix 
	and $\H^2$ denotes the hyperbolic plane which is the upper half plane $\H^2=\{(x, y)\in\R^2: y>0\}$, endowed  
	with the hyperbolic metric $ds^2=\frac{dx^2+dy^2}{y^2}$. The hyperbolic plane 
	$\H^2$ has constant curvature $-1$. 
	
	 The group $\PSL(2,\R)$  acts transitively on $\H^2$ by 
	M\"obius transformations
	$z\mapsto \frac{az+b}{cz+d}$.
	If the action has no fixed points, then the factor $\Gamma\backslash\H^2$  is a closed Riemann surface of genus at least $2$ 
	and has the hyperbolic plane $\H^2$ as the universal covering.  Furthermore, 
	the natural projection $\pi_\Gamma: \H^2\rightarrow \Gamma\backslash\H ^2, \pi_\Gamma(z)=\Gamma z,\ z\in\H^2$ becomes a local isometry. This implies that $\Gamma\backslash\H^2$ also has constant curvature $-1$. 
	The natural Riemannian metric on $\PSL(2,\R)$ and the  Sasaki metric (see \cite{pater}) on the 
	unit tangent bundle $T^1\H^2$ induce an isometry from $T^1\H^2$ to $\PSL(2,\R)$.
	As a consequence, with the Sasaki metric  on the unit tangent bundle $T^1(\Gamma\backslash\H^2)$ 
	induced by the natural Riemannian metric on $\Gamma\backslash\H^2$ 
	and the natural Riemannian metric on $\Gamma\backslash\PSL(2,\R)=\{\Gamma g,g\in\PSL(2,\R)\}$ that is the collection of the right co-sets of $\Gamma$ in $\PSL(2,\R)$, 
	 there exists an isometry $\Xi$ from the unit tangent bundle $T^1(\Gamma\backslash\H ^2)$
	to the quotient space 
	$\Gamma\backslash \PSL(2,\R)$. Instead of considering $T^1(\Gamma\backslash\H^2)$ it is often easier to work on $\Gamma\backslash {\rm PSL}(2, \R)$, which is a three-dimensional real connected analytic manifold.

   In the hyperbolic plane, a horocycle is a euclidean circle tangent to the real axis or a horizontal line. 
   Let $b_t\in\PSL(2,\R)$ denote the equivalence class obtained from the matrix $B_t=\scriptsize\big(\begin{array}{cc}
   1 & t\\ 0 & 1
   \end{array}\big)\in\SL(2,\R)$. 
   The (stable) horocycle flow on $T^1\H^2$ can be described by the flow $\theta_t^\G(g)=gb_t$ on $\G=\PSL(2,\R)$ (see \cite[Chapter 9]{einsward}).
   A flow $(\theta^X_t)_{t\in\R}$ on $X=\Gamma\backslash\PSL(2,\R)$ defined by
   \begin{equation*} 
   \theta^X_t(x)=\Gamma g b_t,\quad t\in\R, \quad x=\Gamma g\in X
   \end{equation*}    
   is smooth and satisfies $\Pi_\Gamma\circ\theta_t^\G
   =\theta_t^X\circ\Pi_\Gamma$, where $\Pi_\Gamma: \G\to X, \Pi_\Gamma(g)=\Gamma g, g\in\G$ denotes the natural projection. This shows how the flow $(\theta_t^\G)_{t\in\R}$ on $\G$ induces  the `quotient flow' $(\theta_t^X)_{t\in\R}$ on $X$
    and the horocycle flow $(\theta_t^\X)_{t\in\R}$  on $\X=T^1(\Gamma\backslash\H^2)$ can be equivalently described as the flow $(\theta^X_t)_{t\in\R}$
	on $X=\Gamma\backslash\PSL(2,\R)$ 
	by the relation 
	\begin{equation}\label{cr}  
	\theta_t^\X=\Xi^{-1}\circ\theta_t^X\circ\Xi \quad \mbox{for all}\quad t\in\R.
	\end{equation}
	The flow $(\theta^X_t)_{t\in\R}$ is ergodic (see \cite{kai}). In the case that $X$ is compact, this flow
has no periodic points as well as fixed points (see \cite{huynh}), and every its orbit is dense.
	
	For more details of all statements above, we refer the readers to \cite{bedkeanser,einsward}. 
		In order to simplify notation, we will drop the superscript $X$ from $(\theta^X_t)_{t\in\R}$. In the whole present paper, we always assume the action of $\Gamma$ on $\H^2$ to be free (of fixed points) and the factor   $\Gamma\backslash\H^2$  to be compact.  Note that the compactness of  $\Gamma\backslash\H^2$, $T^1(\Gamma\backslash\H^2)$, and $\Gamma\backslash \PSL(2,\R)$ are equivalent.

	In the rest of this section  we collect some notions and auxiliary technical results which will be used
	 afterwards. 

	For $g=[G]\in\PSL(2,\R), G=\big(\scriptsize\begin{array}{cc}a&b\\c&d\end{array} \big)\in\SL(2,\R)$,
the trace of $g$ is defined by \[\tr(g)=|a+d|.\]
If the action of $\Gamma$ on $\H^2$ is free and the factor $\Gamma\backslash\H^2$ is compact
then all elements $g\in\Gamma\setminus\{e\}$ are hyperbolic \cite[Theorem 6.6.6]{ratcliff}, i.e. $\tr(g)>2$;
here $e=[E_2]$ denotes the unity of $\PSL(2,\R)$.  
Furthermore, one gets a stronger result:
\begin{lemma}\label{hy}
	If the factor $\Gamma\backslash\H^2$ is compact, then there exists $\eps_*>0$ such that
	\[\tr(g)\geq 2+\eps_*\quad \mbox{for all}\quad g\in \Gamma\setminus \{e\}. \]
\end{lemma}

\begin{lemma}\label{at}
			
	(a)	There is a natural Riemannian metric on $\G=\PSL(2,\R)$
		such that the induced metric function $d_\G$ is left-invariant under $\G$ and 
		\[d_\G(b_t,e)\leq |t|\quad\mbox{for all}\quad t\in\R. \]
		
	(b)
	For every $\eps>0$ there is $\delta>0$ with the following property. 
	If $g\in\G$ satisfies $d_\G(g,e)<\delta$ then there is $G=\Big(\scriptsize{\begin{array}{cc}g_{11}&g_{12}\\
	g_{21}&g_{22}
	\end{array}}\Big)\in\SL(2,\R)$
	such that $g=[G]$ and
	\[|g_{11}-1|+|g_{12}|+|g_{21}|+|g_{22}-1|<\eps.\]
\end{lemma}

See \cite[Subsection 9.3]{einsward} for more details of (a) and \cite[Lemma 3.4 (b)]{huynh} for a proof of (b).

	We define a metric function $d_{X}$ on $X=\Gamma\backslash\PSL(2,\R)$ by 
	\[ d_{X}(x_1, x_2)
	=\inf_{\gamma_1, \gamma_2\in\Gamma} d_{\G}(\gamma_1 g_1, \gamma_2 g_2)
	=\inf_{\gamma\in\Gamma} d_{\G}(g_1, \gamma g_2), \]   
	where $x_1=\Gamma g_1 $, $x_2=\Gamma g_2$.
	In fact, in the case that $X$ is compact, one can prove that the infimum is a minimum, i.e., 
	for $x_1=\Gamma g_1, x_2=\Gamma g_2\in X$, there exists a $\gamma_0\in\Gamma$ such that
	\begin{equation*}
	d_{X}(x_1, x_2)=\min_{\gamma\in\Gamma} d_{\G}(g_1, \gamma g_2)=d_{\G}(g_1, \gamma_0 g_2).
	\end{equation*} 
		It is possible to derive 
	a uniform lower bound on $d_{\G}(g, \gamma g)$ 
	for $g\in \PSL(2,\R)$ and $\gamma\in\Gamma\setminus\{e\}$ as follows.
	
	\begin{lemma}\label{sigma_0}
		If the space $X=\Gamma\backslash\PSL(2,\R)$ is compact, then there exists $\sigma_0>0$ such that
		\begin{equation}\label{sigma0}d_\G(\gamma g, g)>\sigma_0\quad\mbox{for all}\quad \gamma\in\Gamma\setminus\{e\}. 
		\end{equation}
	\end{lemma}
	See \cite[Lemma 1, p. 237]{ratcliff} for a similar result on $\Gamma\backslash\H^2$.

\section{Concepts and main results}\label{sec3}
In this section we start by recalling the notion of time changes for flows and the definitions of some expansive properties and next we  prove that the horocycle flow $(\theta^\X_t)_{t\in\R}$ on $\X=T^1(\Gamma\backslash\H^2)$ is positive strong kinematic expansive. Finally 
we provide a solution to Flinn's conjecture.

\begin{definition}[\cite{flinn}]\label{tcdn}\rm Let $M$ be a metric space and $\phi,\psi:\R\times M\to M$ be continuous flows. We say 
	that $(\phi_t)_{t\in\R}$ is a {\em time change} of $(\psi_t)_{t\in\R}$ if 
	for every $x\in M$ the orbits $\phi_\R(x)$, $\psi_\R(x)$ 
	and their orientations coincide.
\end{definition}

One has an expression for a time change of a given flow via its respective flow as follows.
\begin{lemma}\label{alpha}Let $M$ be a metric space.
	If $(\phi_t)_{t\in\R}$ is a time change of $(\psi_t)_{t\in\R}$ on $M$, then there is
	a unique 
	continuous function $\alpha: \R\times M\to M$ such that 
	\[\psi_t(x)=\phi_{\alpha(t,x)}(x)\quad\mbox{for all}\quad (t,x)\in \R\times M. \]
\end{lemma}
For a proof, see \cite[Lemma 3.11]{flinn}.

Next, we recall the definitions of expansive (in the sense of Flinn), separating (in the sense of Gura) and kinematic expansive (in the sense of Artigue) flows. In what follows, let $\phi:\R\times M\to M$ be a continuous flow on a compact metric space $(M,d)$.
\begin{definition}[\cite{flinn}] \rm We say that $(\phi_t)_{t\in\R}$ is  {\em expansive}	if for each $\eps>0$, 
	there exists $\delta>0$ with the following property. 
	If $s:\R\rightarrow \R$ is a homeomorphism with $s(0)=0, s(t)>0$ for $ t>0$ and 
	\[d(\phi_t(x),\phi_{s(t)}(y))<\delta \quad \mbox{for all}\quad t\in \R\]
	then $y=\phi_\tau(x)$ for some $\tau\in (-\eps,\eps)$.
\end{definition} 
 In the case that the homeomorphism $s$ is replaced by the identity map, the expansiveness becomes the kinematic expansiveness:
	\begin{definition}[\cite{artigue}]\label{dnke}\rm We say that $(\phi_t)_{t\in\R}$ is {\em kinematic expansive} 
		 if for each $\eps>0$, there exists $\delta>0$ with the following property. If
		\begin{equation}\label{exdn}d(\phi_t(x),\phi_{t}(y))<\delta \quad \mbox{for all}\quad t\in (-\infty,+\infty)	
		\end{equation}
		then $y=\phi_\tau(x)$ for some $\tau\in(-\eps,\eps)$.
	\end{definition}
	If `$t\in (-\infty,+\infty)$' in \eqref{exdn} is replaced by `$t\in [0,+\infty)$'
	then the flow is said to be {\em positive kinematic expansive} (\cite{artigue3}). If the condition  `$\tau\in(-\eps,\eps)$' is ignored, the flow is called positive separating (called {\em unilaterally separating} in \cite{gura}):
	\begin{definition}[\cite{artigue}]\rm We say that $(\phi_t)_{t\in\R}$ is {\em positive separating}
	 if  there exists $\delta>0$ with the following property. If
		\begin{equation}\label{sdn}
		d(\phi_t(x),\phi_t(y))<\delta \quad \mbox{for all}\quad t\in [0,+\infty)
			\end{equation} 
		then $y=\phi_\tau(x)$ for some $\tau\in\R$, i.e. the points $x$ and $y$ lie on the same orbit. 
	\end{definition}
Such a $\delta$ is called a {\em positive separating constant}.
\begin{definition}\rm We say that $(\phi_t)_{t\in\R}$ is 
	{\em strong kinematic expansive} or  {\em weakly expansive}
	if every time change of $(\phi_t)_{t\in\R}$ is  kinematic expansive.
	
	{\em Positive strong expansive} and {\em strong separating expansive} flows are defined analogously.
\end{definition}

\begin{definition}[\cite{artigue}]\rm Let $X,Y$ be metric spaces.
	Two continuous flows $\phi:\R\times X\to X$ and $\psi: \R\times Y\to Y$ are said to be 
	{\em equivalent} if there is a homeomorphism $h:X\to Y$ such that $\phi_t=h^{-1}\psi_th$ for all $t\in\R$. 
\end{definition}
Due to \eqref{cr}, the flows $(\theta^\X_t)_{t\in\R}$ and $(\theta_t)_{t\in\R}$ are equivalent.  
It is easy to see that all the expansive properties introduced above are invariant under equivalence.

Next, we recall the positive strong separation for the horocycle flow showed by Gura \cite{gura}. 
\begin{lemma}[\cite{gura}]\label{ss} The horocyle flow on a compact oriented connected surface
	provided with a Riemannian metric  of negative curvature of class $C^3$ is positive strong separating.
\end{lemma}

Now we are in a position to state the first main result of this paper. The following gives a stronger property of the horocyle flow  on (the unit tangent bundle of) $\Gamma\backslash\H^2$ which is a compact oriented connected analytic Riemann surface with constant negative curvature. 
	
	\begin{theorem}\label{pne}
		The horocycle flow $(\theta_t^\X)_{t\in\R}$ on $\X=T^1(\Gamma\backslash\H^2)$ is positive strong kinematic expansive.
	\end{theorem}
	{\bf Proof\,:}
Since the positive strong kinematic expansiveness is invariant under equivalence, 
 it suffices to show that the flow $(\theta_t)_{t\in\R}$ is positive strong kinematic expansive.
	Let $(\psi_t)_{t\in\R}$ be a time change of the flow $(\theta_t)_{t\in\R}$ (see Definition \ref{tcdn}). 
	By Lemma \ref{alpha},  there exist  continuous functions $\alpha,\beta: \R\times X\to\R $ satisfying 
	\[\theta_t(x)=\psi_{\alpha(t,x)}(x)\quad\mbox{and}\quad \psi_t(x)=\theta_{\beta(t,x)}(x) \quad \mbox{for all}\quad (t, x)\in\R\times X. \]
	For every $\eps>0$, there is $\rho=\rho(\eps)\in (0,\eps_*)$ so that
	\begin{equation}\label{stt2} |\alpha(t,x)|<\eps\quad\mbox{for all}\quad (t,x)\in (-\rho,\rho)\times X;
	\end{equation} 
	recall $\eps_*$ from Lemma \ref{hy}. 
	For $\rho>0$ above, take $\delta=\delta(\rho)\in(0,\min\{\sigma_0/4,\eps_*\})$ as in Lemma \ref{at}\,(b),
	with $\sigma_0$ from Lemma \ref{sigma_0}. According to Lemma \ref{ss}, the flow $(\theta_t)_{t\in\R}$ is positive strong separating. This implies that the time change 
	$(\psi_t)_{t\in\R}$ is positive separating.
	Let $0<\varrho<\delta$ be a positive separating constant. For any $x,y\in X$, if 
	\begin{equation}\label{stt3}
	d_X(\psi_t(x),\psi_t(y))<\varrho\quad
	\mbox{for} \quad t\geq 0
	\end{equation}
	then $y=\psi_r(x)$ for some $r\in \R$. We need to show that $|r|<\eps$. 
		There exists a unique $s\in \R$ such that $r=\alpha(s,x)$ and  $y=\theta_s(x)$.
	Due to \eqref{stt2}, it remains to verify  $|s|<\rho$.
	
	First, we rewrite \eqref{stt3} as 
	\[ d_X(\theta_{\beta(t,x)}(x),\theta_{\beta(t,y)+s}(x))<\varrho\quad\mbox{for}\quad t\geq 0. \]
	Let $s_1(t)=\beta(t,x), s_2(t)=\beta(t,y)+s, t\in \R$ and write $x=\Gamma g$ for $g\in\G$. 
	For every $t\geq 0$, there is $\gamma(t)\in\Gamma$ so that
	\begin{equation*}\label{strong1}
	d_X(\theta_{s_1(t)}(x),\theta_{s_2(t)}(x))=d_\G(\gamma(t) gb_{s_1(t)}, gb_{s_2(t)})<\varrho. 
	\end{equation*} 
	We claim that  $\gamma(t)=\gamma(0)=:\gamma$ for all $t\geq 0$. Indeed, for any $t_1,t_2\in[0,+\infty)$, we have
	\begin{eqnarray*}
		\lefteqn{d_\G(\gamma(t_2)^{-1}\gamma(t_1)gb_{s_1(t_1)},gb_{s_1(t_1)})}
		\\
		&=& d_\G(\gamma(t_1)gb_{s_1(t_1)},\gamma(t_2)gb_{s_1(t_1)})
		\\
		&\leq& d_\G(\gamma(t_1)gb_{s_1(t_1)},gb_{s_2(t_1)})
		+d_\G(gb_{s_2(t_1)},gb_{s_2(t_2)})\\
		&&+d_\G(gb_{s_2(t_2)},\gamma(t_2)gb_{s_2(t_2)})
			+\ d_\G(\gamma(t_2)gb_{s_2(t_2)},\gamma(t_2)gb_{s_1(t_1)})
		\\
		&=& d_\G(\gamma(t_1)gb_{s_1(t_1)},gb_{s_2(t_1)})
		+d_\G(b_{s_2(t_1)},b_{s_2(t_2)})
		\\
		&&
		+d_\G(gb_{s_2(t_2)},\gamma(t_2)gb_{s_1(t_2)})
		+d_\G(b_{s_1(t_2)},b_{s_1(t_1)} )
		\\
		&\leq& 2\varrho+|s_2(t_1)-s_2(t_2)|+|s_1(t_1)-s_1(t_2)|
	\end{eqnarray*}
	due to Lemma \ref{at} (a).
	For a given $L>0$, we verify that $\gamma(t)=\gamma(0)$ for all $t\in [-L,L]$. Since $s_1,s_2:[0,L]\to\R$ is uniformly continuous, there is $0<\eta=\eta(L,\delta)<\delta$
	such that if $t_1,t_2\in [0,L]$ and $|t_1-t_2|\leq\eta$ then 
	$|s_i(t_1)-s_i(t_2)|<\delta$ for $i=1,2$.
	For $t_1,t_2\in [0,\eta]$, then $|t_1-t_2|\leq\eta$ implies $|s_i(t_1)-s_i(t_2)|<\delta$ for $i=1,2$. 
	This yields \[d_\G(\gamma(t_2)^{-1}\gamma(t_1)gb_{s_1(t_1)},gb_{s_1(t_1)})<4\delta<\sigma_0. \]
	From the property of $\sigma_0$ (see \eqref{sigma0}), it follows that 
	$\gamma(t_2)=\gamma(0)$ for all $t_2\in [0,\eta]$. Then we repeat the argument for
	$t_1=\eta$ and $t_2\in [\eta, 2\eta]$, we deduce that $\gamma(t)=\gamma(0)$ for all $t\in [0,2\eta]$, which upon further iteration leads to $\gamma(t)=\gamma(0)$ for all $t\in [0,L]$. Since $L$ is arbitrary, we have
	$\gamma(t)=\gamma(0)$ for all $t\in [0,+\infty)$. Therefore
  	\[d_\G(\gamma gb_{s_1(t)}, gb_{s_2(t)})<\varrho<\delta\quad \mbox{for}\quad t\geq 0. \]
	or
	\[d_\G(b_{-s_2(t)}g^{-1}\gamma gb_{s_1(t)}, e)<\delta\quad \mbox{for}\quad t\geq 0. \]
	According to Lemma \ref{at}\,(b), there is $K=\big(\scriptsize\begin{array}{cc} k_{11}&k_{12}\\k_{21}&k_{22}\end{array}\big)\in\SL(2,\R)$ such that $g^{-1}\gamma g=[K]$ and
	\begin{align}\notag
		|k_{11}-k_{21}s_2(t)-1 |
		&+ | (k_{11}-k_{21}s_2(t))s_1(t)-k_{22}s_2(t)+k_{12}| \\
		&
	+ | k_{21}|+ |k_{21}s_1(t)+k_{22}-1|<\rho\quad\mbox{for}\quad t\geq 0,  \label{st1}
	\end{align}
	using
	\begin{eqnarray*}
	 B_{-s_2(t)}KB_{s_1(t)}=\Big(\begin{array}{cc} k_{11}-k_{21}s_2(t)&(k_{11}-k_{21}s_2(t))s_1(t)-k_{22}s_2(t)+k_{12}\\
	k_{21}&k_{21}s_1(t)+k_{22}\end{array}\Big).
	\end{eqnarray*}
	Noting that $s_1(t)\to +\infty$ and $s_2(t)\to+\infty$ as $t\to +\infty$, 
	it follows from \eqref{st1}  that $k_{21}=0$ and
	\begin{equation}\label{st2} |k_{11}-1|+|k_{11}s_1(t)-k_{22}s_2(t)+k_{12}|+|k_{22}-1|<\rho\quad\mbox{for}\quad t\geq 0. 
	\end{equation} 
	This yields 
	$\tr(\gamma)=\tr(g^{-1}\gamma g)=|k_{11}+k_{22}|<2+\rho<2+\eps_*$. 
	Owing to Lemma \ref{hy}, we have $\gamma=e$ and hence  $k_{11}=k_{22}=1, k_{12}=k_{21}=0$.
	Then \eqref{st2} implies $|s_1(t)-s_2(t)|<\rho$ for all $t\geq 0$. Finally, $s_1(0)=0$ and $s_2(0)=s$ yield $|s|<\rho$, which completes the proof.
	{\hfill$\Box$}

\bigskip

The remainder of this paper is devoted to solving Flinn's conjecture. In his thesis,
Flinn  showed that the expansiveness implies the weak expansiveness  (see \cite[Theorem 3.24]{flinn}). In addition, he conjectured that the converse is also true:

\noindent
{\bf Flinn's conjecture:} {\em  Let  $\phi:\R\times M\to M$
be a continuous flow on a compact metric space $M$. 
If $(\phi_t)_{t\in\R}$ is weakly expansive then it is also expansive. 
}

\begin{theorem}\label{fca} Flinn's conjecture is not true.
\end{theorem}
{\bf Proof\,:} We show that the horocycle flow $(\theta_t)_{t\in\R}$ is not expansive, although
	it is weakly expansive according to Theorem \ref{pne}. To this end, we recall Remark 2.8 (a) in \cite{huynh}. For any $\delta>0$, 
	take $\rho=\rho(\delta)$ as in Lemma \ref{at}\,(b).
	Let $a>0$ be such that $|a-1|+|1/a-1|<\rho$ and  $|a+1/a|< 2+\eps_*$ with $\eps_*$ from Lemma \ref{hy}.
	Let $x=\Gamma e, y=\Gamma h$ with $h=[H], H=\Big(\begin{array}{cc}a&0\\ 0& 1/a\end{array}\Big)$.
	Set $s(t)=a^2 t$,  $t\in\R$ to have
	\begin{eqnarray*}
		d_X(\theta_{s(t)}(x),\theta_t(y))
		&=&d_X(\Gamma b_{s(t)},\Gamma hb_t)\leq d_\G(b_{s(t)},hb_t )=d_\G(b_{-s(t)}hb_{t},e)\\
		&=&d_\G(h,e)<\delta\quad\mbox{for all}\quad t\in\R.
	\end{eqnarray*}
		But $x$ and $y$ are not in the same orbit since otherwise there would exist  $\tau\in\R$ such that $y=\theta_\tau(x)$. Then $\gamma = h b_\tau$ for some $\gamma\in\Gamma$ implies $\tr(\gamma)=|a+1/a|<2+\eps_*$.
		It follows from Lemma \ref{hy} that $\gamma=e$ and hence $b_{-\tau}=h=[H]$, which is impossible. 
		We have showed that  the horocycle flow is not expansive in the sense of Flinn, although it is strong kinematic expansive,
	i.e., weakly expansive. This gives a negative answer to Flinn's conjecture.
{\hfill$\Box$}
	\bigskip 
	
\noindent	{\bf Acknowledgments:} This work was supported by the Research Project of
Vietnam Ministry of Education and Training (Grant No. B2019-DQN-10). I would like to thank an anonymous referee for carefully reading my paper and suggesting me solve Flinn's conjecture. I enjoyed many fruitful discussions with Alfonso Artigue.

\end{document}